\documentclass[preprint,sort&compress,final,12pt]{elsarticle}

\usepackage{amsmath}
\usepackage{amssymb}
\usepackage{graphicx}
\usepackage{latexsym}
\usepackage{epstopdf}
\usepackage{natbib}
\usepackage{color}
\usepackage{hyperref}

\newtheorem{theorem}{Theorem}[section]

\newtheorem{lemma}[theorem]{Lemma}

\numberwithin{equation}{section}

\def\Proof{\noindent{\bf Proof.}~}
\def\qed{\hfill$\square$\smallskip}

\def\dsum{\displaystyle\sum}

\def\Im{\mathrm{Im}}

\journal{\empty}
\date{}

\begin{document}

\begin{frontmatter}

\title{The linearization of periodic Hamiltonian systems with one degree of freedom under the Diophantine condition}

\author[au1,au2]{Nina Xue}

\address[au1]{School of Mathematical Sciences, Beijing Normal University, Beijing 100875, P.R. China.}

\author[au1]{Xiong Li\footnote{ Partially supported by the NSFC (11571041) and the Fundamental Research Funds for the Central Universities. Corresponding author.}}

\address[au2]{School of Mathematics and Information Sciences, Weifang University, Weifang, Shandong, 261061, P.R. China.}

\ead[au2]{xli@bnu.edu.cn}

\begin{abstract}
In this paper we are concerned with the periodic Hamiltonian system with one degree of freedom, where the origin is a trivial solution. We assume that the corresponding linearized system at the origin is elliptic, and the characteristic exponents of the linearized system are $\pm i\omega$ with $\omega$ be a Diophantine number, moreover if the system is formally linearizable, then it is analytically linearizable. As a result, the origin is always stable in the sense of Liapunov in this case.
\end{abstract}


\end{frontmatter}

\section{Introduction}
We consider a Hamiltonian system with one degree of freedom
\begin{equation} \label{intro1}
\dot{x}=\frac{\partial H}{\partial y}(x,y,t),\ \ \ \ \dot{y}=-\frac{\partial H}{\partial x}(x,y,t),
\end{equation}
where $H(0,0,t)=H_x(0,0,t)=H_y(0,0,t)=0$, the dot indicates differentiation with respect to the time $t$. The Hamiltonian function $H$ is continuous and $2\pi$-periodic in $t$, real analytic in a neighborhood of the origin $(x,y)=(0,0)$, and the Taylor series of $H$ in a neighborhood of the origin is assumed to be
\begin{equation} \label{intro2}
H(x,y,t)=H_{2}+H_{3}+\cdots +H_{j}+ \cdots,
\end{equation}
where
$$H_{j}=\sum_{\mu+\nu=j}h_{\mu\nu}(t)x^{\mu}y^{\nu}, \ \ \ \ j=2,3,\cdots,$$
and the coefficients $h_{\mu\nu}(t)$ are $2\pi$-periodic with respect to the time $t$.

Initially we want to investigate the problem about the stability of the trivial solution of the periodic Hamiltonian system (\ref{intro1}). There are plenty of works about the stability of the trivial solution, and we can refer to \cite{Moser1}, \cite{Siegel} for a detailed description. For recent developments, one may consult \cite{Liu1}, \cite{Ortega3} and the references therein. For time-periodic Lagrangian equations, an analytical method called the third order approximation, has been developed recently by Ortega in a series of papers \cite{Ortega0}, \cite{Ortega1}, \cite{Ortega2}, \cite{Ortega4}. After that, some researcher have extend the applications of the third approximation, and some stability results for several types of Lagrangian equations have been established. We refer the reader to \cite{Lei}, \cite{Ortega3} for the forced pendulum equation, and \cite{Lei}, \cite{Torres1}, \cite{Torres2} for some singular equations.

It is well known that for periodic Hamiltonian system (\ref{intro1}), if there exists a characteristic exponent of the linearized system with non zero real part, then the trivial solution is unstable, thus we assume that the linearized system is stable and that the characteristic exponents are pure imaginary, say $\pm i\omega.$

In the general elliptic case, Arnold \cite{Arnold} and Moser \cite{Moser}, \cite{Siegel} proved that the solution is stable if certain nondegeneracy conditions are fulfilled. In the case of resonance, i.e., $\omega$ is a rational number, the above result cannot be applied directly. For this case, there also are many results, see \cite{Bardin}, \cite{Caraballo}, \cite{Mansilla}, \cite{Markeyev}, \cite{Santos1}, \cite{Santos2}, \cite{Santos3} and the references therein. For example, Mansilla \cite{Mansilla} obtained some sufficient conditions for stability and instability of the trivial solution by using Moser's twist theorem and Liapunov theorem, respectively. Bardin \cite{Bardin} studied the degenerate case by using Lie normal form, and established some general criteria to solve the stability problem. The basic method of these results is that by considering up to terms of a certain order in the normal form, then the sufficient conditions for stability and instability are obtained, respectively.

On the other hand, R\"{u}ssmann \cite{Russmann1} considered a real analytic area-preserving mapping near the fixed point $(0,0)$ of the form
\begin{equation}\label{R1}
\left\{
\begin{array}{ll}
x_{1}=f(x,y)=x\cos\gamma_{0}-y\sin\gamma_{0}+\cdots,\\[0.2cm]
y_{1}=f(x,y)=x\sin\gamma_{0}+y\cos\gamma_{0}+\cdots,
\end{array}
\right.
\end{equation}
where $f$ and $g$ are convergent power series in $x, y$ with real coefficients and $\frac{\gamma_{0}}{2\pi}$ is an irrational real number. He proved that if \eqref{R1} is formally linearizable, then it is analytically linearizable under the Bruno condition. We recall the main ideas of the proof. R\"{u}ssmann first studies the process of formal normalization using a functional iterative approach, then he constructs a formal iteration process converging to a zero of the operator $\mathcal{F}:=f\circ\varphi-\varphi\circ\Lambda$ (where $\Lambda$ is the linear part of the area-preserving mapping), and assuming $\Lambda$ diagonal, he gives estimates for each iteration step, proving that, under the Bruno condition, the process converges to a holomorphic linearization. Recently, J. Raissy \cite{Raissy2} had given a new proof of the above R\"{u}ssmann's result under a slightly different arithmetic hypothesis, and using direct computation via power series expansions in presence of resonances.

Motivated by the above results, especially by R\"{u}ssmann's work, we will consider the linearization of the periodic Hamiltonian system \eqref{intro1}. Let us assume that $\omega$ is an irrational real number. Without loss of generality we assume (see \cite{Moser0} for more details), that a linear $2\pi$-periodic symplectic transformation has already been made such that $H_{2}$ is given by
$$H_{2}=\frac{\omega}{2}(x^{2}+y^{2}).$$
Introduce symplectic action angle variables by means of the formula
$$x=\sqrt{2r}\cos\theta, \ y=\sqrt{2r}\sin\theta.$$
In the variables $r$ and $\theta$ , the Hamiltonian \eqref{intro2} reads as
\begin{equation} \label{intro3}
H(r,\theta,t)=H_{2}+H_{3}+\cdots +H_{j}+ \cdots,
\end{equation}
where
$$H_{2}=\omega r, \
H_{j}=\sum_{\mu+\nu=j}h_{\mu\nu}(t)(2r)^\frac{j}{2}\cos^{\mu}\theta \sin^{\nu}\theta, \ \ j\geq 3,$$
and $\omega$ is an irrational real number.

By Lie process of normalization, there exists a formal, symplectic, $2\pi$-periodic change of variables which transforms \eqref{intro3} into the new Hamiltonian
\begin{equation*}
\widetilde{H}(\tilde{r},\tilde{\theta},t)=\omega \tilde{r}+A_{4}\tilde{r}^{2}+A_{6}\tilde{r}^{3}+\cdots +A_{2j}\tilde{r}^{j}+ \cdots,
\end{equation*}
where $A_{2j}\, ( j=2,3,\cdots )$ are constants, and the details can be found in Section 2 of this paper. For the sake of brevity, we omit the tildes on the letters, and the new Hamiltonian has the form
\begin{equation} \label{intro4}
H(r,\theta,t)=\omega r+A_{4}r^{2}+A_{6}r^{3}+\cdots +A_{2j}r^{j}+ \cdots.
\end{equation}

If there exists  $A_{2j}\neq 0 $ for some $ j=2,3,\cdots ,$ it follows by Moser's twist theorem, the origin in \eqref{intro1} is stable in the sense of Liapunov. If all constants $A_{2j}=0\, ( j=2,3,\cdots ),$ then we have the following theorem.

\begin {theorem}\label{main th}
If there is a formal, symplectic, $2\pi$-periodic change of variables which transforms the given Hamiltonian \eqref{intro3} into the linearized system
\begin{equation} \label{intro5}
H=H_{2}=\omega r,
\end{equation}
and there exist constants $\alpha >0, \ \tau>1 $ such that
\begin{equation} \label{intro6}
|k\omega -l|\geq \frac{\alpha}{|k|^{\tau}}
\end{equation}
for all integers $k,l$ with $k\neq 0$, then there is also an analytic, symplectic, $2\pi$-periodic change of variables which transforms \eqref{intro3} into \eqref{intro5}.

\end{theorem}

A consequence of this theorem is that under the condition \eqref{intro6}, in the case that all constants $A_{2j}=0 \ (j=2,3,\cdots)$,  the origin in \eqref{intro1} is also stable. Therefore if the frequency $\omega$ in the linearized system \eqref{intro5} satisfies the Diophantine condition \eqref{intro6}, the origin in \eqref{intro1} is always stable.

Other points of view or approaches related to the direction of linearization can be found in different articles. For example, J. Raissy \cite{Raissy1} considered the linearization of holomorphic germs with quasi-Brjunu fixed points, and proved that, under certain arithmetic conditions on the eigenvalues of $Df_O$ (where $f$ is a germ of holomorphic diffeomorphism of $\mathbb{C}^n$ fixing the origin $O$) and some restrictions on the resonances, $f$ is locally holomorphically linearization if and only if there exists a $f$-invariant complex manifold. J. Raissy in \cite{Raissy3} investigated the shape that a simultaneous linearization of the given germs can have, and proved that if $f_1, \cdots, f_n$ (where $f_1, \cdots, f_n$ are $h\geq 2$ germs of biholomorphisms of $\mathbb{C}^n$ fixing the origin) commute and their linear parts are almost simultaneously Jordanizable, then they are simultaneously formally linearizable. He also proved that if the linear terms of the germs are diagonalizable, then they are holomorphically simultaneously linearizable under Brjuno-type condition. Wu in \cite{Wu} studied the linearizability of a class quasi-periodic Hamiltonian systems. He proved that, under some conditions, if the system is formally linearizable, then it is analytically linearizable in the resonant case.

This paper is organized as follows. In Section 2, we give some basic results about Lie normal form that will be useful. In Section 3, we describe the idea of the proof. Then we give the complete proof of the main result in Section 4 and Section 5.

\section{Lie normal form}
In this section, we introduce the Lie normal form theory for a general Hamiltonian system of $n$ degrees of freedom. According to the well-known concept of normal form, for a given Hamiltonian system, we try to change coordinate so that the system of the new coordinate is simpler. The definition of the simplicity depends on the problem at hands. For Hamiltonian systems, we often use the method of Lie transforms developed by Deprit which is a general procedure to change variables in a Hamiltonian system that depends on a small parameter.

Consider the system
\begin{equation} \label{nm1}
\dot{x}=J\nabla H(\varepsilon,t,x),
\end{equation}
where
$x^{T}=(x_{1}, \cdots, x_{2n}),
\nabla H^{T}=\left(\frac{\partial H}{\partial x_{1}},\cdots, \frac{\partial H}{\partial x_{2n}}\right),$
and
\begin{equation*}
J=\left(
\begin{matrix}
0&I&\\
-I&0&
\end{matrix}
\right),
\end{equation*}
where $0$ is the $n\times n$ zero matrix and $I$ is the $n\times n$ identity matrix.
Assume that $H$ has series expansion in the small parameter $\varepsilon$, that is,
\begin{equation} \label{nm2}
H(\varepsilon,t,x)=H_{\ast}(\varepsilon,t,x)=
\sum_{i=0}^{\infty}\left(\frac{\varepsilon^{i}}{i!}\right)H_{i}^{0}(t,x),
\end{equation}
where  $H^{0}_{i}$ is a homogeneous polynomial in $x$ of degree $i+2$.
Make a symplectic change of coordinates, $x=X(\varepsilon,t,y),$ which is the solution of the Hamiltonian system
$$\frac{dx}{d\varepsilon}=J\nabla W(\varepsilon,t,x)$$
satisfying the initial condition
$$x(0)=y.$$
The change of coordinates transforms \eqref{nm1} into the Hamiltonian system
$$\dot{y}=J\nabla G(\varepsilon,t,y)+J\nabla \widetilde{G}(\varepsilon,t,y):=J\nabla K(\varepsilon,t,y),$$
where $G(\varepsilon,t,y)=H(\varepsilon,t,X(\varepsilon,t,y))$ is the Lie transform of $H$, and
$$\widetilde{G}=-\int_{0}^{\varepsilon}\frac{\partial W}{\partial t}(s,t,X(s,t,y))ds$$
is the remainder function, and $K:=G+\widetilde{G}$ is the new Hamiltonian.
Let $K,W$ have series expansions of the form
\begin{equation} \label{nm3}
K(\varepsilon,t,y)=H^{\ast}(\varepsilon,t,y)=\sum_{i=0}^{\infty}
\left(\frac{\varepsilon^{i}}{i!}\right)H^{i}_{0}(t,y),
\end{equation}

\begin{equation} \label{nm4}
W(\varepsilon,t,x)=\sum_{i=0}^{\infty}
\left(\frac{\varepsilon^{i}}{i!}\right)W_{i+1}(t,x),
\end{equation}
where  $H^{i}_{0}$ is a homogeneous polynomial in $x$ of degree $i+2$, $W_{i}$ is a homogeneous polynomial in $x$ of degree $i$.

According to the Deprit's general algorithm, the following general perturbation conclusion holds.

\begin{lemma}\label{nm1}(See [\cite{Meyer}, Chapter 10].)
Let $\{\mathcal{P}_{i}\}_{i=0}^{\infty}, \ \{\mathcal{Q}_{i}\}_{i=0}^{\infty}$ and $\{\mathcal{R}_{i}\}_{i=0}^{\infty}$ be sequences of linear spaces of smooth functions defined on a common domain $ O $ in $\mathbb{R}^{1}\times\mathbb{R}^{2n} .$ Let $\dot{\mathcal{R}}_{i}$ be the space of all derivatives of functions in $\mathcal{R}_{i}.$ Assume that\\
(1) $ \mathcal{Q}_{i}\subset \mathcal{P}_{i}, \ i=1,2,\cdots;$\\
(2) $H_{i}^{0}\in \mathcal{P}_{i}, \ i=1,2,\cdots;$\\
(3) $ \{\mathcal{P}_{i}, \mathcal{R}_{j}\}\subset \mathcal{P}_{i+j} $ and
$\{\mathcal{P}_{i}, \dot{\mathcal{R}}_{j}\}\subset \mathcal{P}_{i+j},
\ i+j=1,2,\cdots;$
where $\{\cdot ,\cdot \}$ is the Poisson bracket, that is,
$$\{F,G\}=\sum_{i=1}^{n}\left(\frac{\partial F}{\partial x_{i}}\frac{\partial G}{\partial x_{i+n}}
-\frac{\partial F}{\partial x_{i+n}}\frac{\partial G}{\partial x_{i}}\right);$$\\
(4) for any $D\in \mathcal{P}_{i}, \ i=1,2,\cdots,$ there exist $E\in \mathcal{Q}_{i}$ and $C\in \mathcal{R}_{i}$ such that
$$E=D+\{H^{0}_{0},C\}-\frac{\partial C}{\partial t}.$$

Then there exists a $W$ with a formal Hamiltonian of the form  \eqref{nm4}    with $W_{i}\in \mathcal{R}_{i}, i=1,2,\cdots,$ which generates a near-identity symplectic change of variables $x\rightarrow y$ such that the Hamiltonian in the new variables has a series expansion given by  \eqref{nm3}  with $H_{0}^{i}\in \mathcal{Q}_{i}, i=1,2,\cdots.$
\end{lemma}

Now we consider a periodic system
\begin{equation} \label{nm5}
\dot{x}=J\nabla H_{\sharp}(t,x), \ \
H_{\sharp}(t,x)=\sum_{i=0}^{\infty}H_{i}(t,x),
\end{equation}
where $H_{i}$ is a homogeneous polynomial in $x$ of degree $i+2$ with $2\pi$-periodic coefficients. By Floquet theorem, there is a linear, symplectic, periodic change of variables that makes the linear part of the Hamiltonian system be constant in $t$. Therefore we  assume
$H_{0}(t,x)=\frac{1}{2}x^{T}Sx,$ where $S$ is a $2n\times2n$ real constant symmetric matrix, and $A=JS$ is a constant, real, Hamiltonian matrix.

For the periodic Hamiltonian system \eqref{nm5}, in order to study the solutions near the origin, we scale it by $x\rightarrow \varepsilon x$. This is a symplectic transformation with the multiplier $\varepsilon^{-2}$, thus the Hamiltonian becomes
$$H(\varepsilon,t,x)=H_{\ast}(\varepsilon,t,x)=
\sum_{i=0}^{\infty}\left(\frac{\varepsilon^{i}}{i!}\right)H_{i}^{0}(t,x),$$
where $H_{i}^{0}=i!H_{i}.$ Because we are working formally, we can set $\varepsilon=1$ at the end. By Lemma 2.1, we must define spaces $\mathcal{P}_{i},\mathcal{Q}_{i}$ and $\mathcal{R}_{i}$ with $\mathcal{Q}_{i} \subset\mathcal{P}_{i}, H^{0}_{i}\in \mathcal{P}_{i}, H_{0}^{i}\in \mathcal{Q}_{i}.$
The Lie equation to be solved is
$$E=D+\{H^{0}_{0},C\}-\frac{\partial C}{\partial t},$$
where $D$ is given in $\mathcal{P}_{i}$,  we are to find $E\in \mathcal{Q}_{i}$ and
$C\in \mathcal{R}_{i}.$

Suppose $A$ is  simple, let $\mathcal{P}_{i} $ be the space of polynomials in $x$ with coefficients that are smooth $2\pi$-periodic function of $t$. Let $\mathcal{F}=\{H_{0}^{0},\cdot\}-\frac{\partial}{\partial t}:\mathcal{P}_{i} \rightarrow \mathcal{P}_{i}.$ Define $\mathcal{Q}_{i}=\ker(\mathcal{F}),\mathcal{R}_{i}=\mbox{Range}(\mathcal{F})$. Since $A$ is  simple, then $\mathcal{P}_{i}=\mathcal{Q}_{i} \bigoplus \mathcal{R}_{i} .$ Therefore, given $D$ in $\mathcal{P}_{i}$, the Lie equation has a unique $2\pi$-periodic solution $C\in \mathcal{R}_{i}$ and $E\in \mathcal{Q}_{i}$, where $E$ is a $2\pi$-periodic solution of the equation
\begin{equation} \label{nm6}
0=\{H_{0}^{0},E\}-\frac{\partial E}{\partial t}.
\end{equation}
Characterizing the $2\pi$-periodic solution of \eqref{nm6} defines the normal form. Therefore we have the following lemma.

\begin{lemma}\label{nm2}(See [\cite{Meyer}, Chapter 10].)
Suppose  $A$ is  simple, let $H^{0}(x)=H_{0}(x)=\frac{1}{2}x^{T}Sx,$ where $A=JS$ is a constant Hamiltonian matrix. Then there exists a formal, symplectic, $2\pi$-periodic change of variables $x=X(t,y)=y+\cdots $, which transforms the Hamiltonian \eqref{nm5} into the Hamiltonian system
$$\dot{y}=J\nabla H^{\sharp}(t,y), \ \ H^{\sharp}(t,y)=\sum_{i=0}^{\infty}H^{i}(t,y),$$
where
$$\{H_{0},H^{i}\}-\frac{\partial H^{i}}{\partial t}=0, \  \ i=1,2,3\cdots.$$
\end{lemma}

Now we return to the system \eqref{intro1} which is considered in this paper. By Lemma 2.2,
$H_{m}$ \ ($m$ is a fixed integer greater than 2) is in Lie normal form whenever
\begin{equation} \label{nm7}
\{H_{2},H_{m}\}-\frac{\partial H_{m}}{\partial t}=0 .
\end{equation}
Since $H_{2}=\omega r$, it follows from \eqref{nm7} that
\begin{equation} \label{nm8}
\omega \frac{\partial H_{m}}{\partial \theta}=\frac{\partial H_{m}}{\partial t}.
\end{equation}
Assume that
$$ H_{m}(r,\theta,t)=\sum_{k+l=m}h_{kl}(t)e^{i(k-l)\theta}r^{\frac{m}{2}},$$
then equation \eqref{nm8} gives us
$$\sum_{k+l=m} h_{kl}(t)i\omega(k-l)e^{i(k-l)\theta}r^{\frac{m}{2}}
=\sum_{k+l=m}h'_{kl}(t)e^{i(k-l)\theta}r^{\frac{m}{2}}.$$
From this we deduce that if $ h_{kl}(t)\neq 0,$ then $h_{kl}(t)$ satisfies the differential equation
$$ h'_{kl}(t)=i\omega (k-l)h_{kl}(t),$$
whose solutions are
$$h_{kl}(t)=C_{kl}e^{i\omega (k-l)t}, C_{kl}\in \mathbb{C}.$$
Thus we can rewrite $H_{m}$ as
$$H_{m}(r,\theta,t)=\sum_{k+l=m}C_{kl}e^{i(k-l)(\theta+\omega t)}r^{\frac{m}{2}}.$$

Since $H_{m}$ is a real function in its original variables, by symmetry of subindex $k+l=m$, we get
that $C_{kl}=\overline{C}_{lk}$. Moreover, since $H_{m}$ is $2\pi$-periodic in $t$, then
$H_{m}\neq 0$ if and only if $(k-l)\omega\in \mathbb{Z},$ which implies that $k=l$, because $\omega$ is an irrational real number. Then the Hamiltonian in Lie normal form is the form of \eqref{intro4}.

\section{Outline of the proof}

At first we introduce a few notations. Let $D_{\rho,\gamma}=\{r: |r|\leq \rho\}\times \{\theta: |\Im\theta|\leq \gamma \}\subset\mathbb{C}\times\mathbb{C}$, where $|\cdot|$ stands for the sup-norm. The sup-norm of functions on $D_{\rho,\gamma}$ is denoted by $|\cdot|_{\rho,\gamma}.$

The Hamiltonian in \eqref{intro3} will be denoted by
\begin{equation}\label{Originalequ}
H(r,\theta,t)=H_{2}+P,
\end{equation}
where
$$H_{2}=\omega r,\ \ \ \  P=\dsum^{\infty}_{j=3}H_{j}=\sum^{\infty}_{j=3}
\sum_{\mu+\nu=j}h_{\mu\nu}(t)r^{\frac{j}{2}}cos^{\mu}\theta sin^{\nu}\theta,$$
and $P$ is analytic on $D_{\rho,\gamma}$.

By Deprit's method, one can construct a symplectic, $2\pi$-periodic transformation $\varphi$ by means
of the time 1-map of a Hamiltonian flow $X_{F}^{\varepsilon}$, i.e., $\varphi =X_{F}^{1}$, which transforms the Hamiltonian $H$ in \eqref{Originalequ} to
$$H_{+}(r,\theta,t)=G+\tilde{G} , $$
where
$$\ G=H\circ X_{F}^{1}, \ \
 \ \widetilde{G}=-\int^{1}_{0}\frac{\partial F}{\partial t}\circ X_{F}^{\varepsilon}d\varepsilon.$$
Expanding $H_{+}$ in $\varepsilon$ near $\varepsilon=0 $ yields that
\begin{eqnarray*}\label{eq:solu-X2}
 H_{+}&=&H_{2}+\{H_{2},F\}+P+\int^{1}_{0}(1-\varepsilon)\{\{H_{2},F\},F\}\circ X_{F}^{\varepsilon}d\varepsilon  \\
&&{} + \int^{1}_{0}\{P,F\}\circ X_{F}^{\varepsilon}d\varepsilon-\frac{\partial F}{\partial t}-
\int^{1}_{0}\left\{\frac{\partial F}{\partial t},F\right\}\circ X_{F}^{\varepsilon}d\varepsilon  \\
 &=&H_{2}+\{H_{2},F\}+P-\frac{\partial F}{\partial t}\\
 &&{}+\int^{1}_{0}\left\{(1-\varepsilon)\{H_{2},F\}+P-\frac{\partial F}{\partial t},F\right\} \circ X_{F}^{\varepsilon}d\varepsilon.
\end{eqnarray*}
The point is to find $F$ such that
$$ H_{2}+\{H_{2},F\}+P-\frac{\partial F}{\partial t}=H_{2+}$$
is again a normal form.

We want to solve
$$ \{F,H_{2}\}+\widehat{H}_{2}+\frac{\partial F}{\partial t}=P, \ \ \widehat{H}_{2}=H_{2+}-H_{2},$$
namely
\begin{equation}\label{add1}
\{F,H_{2}\}+\frac{\partial F}{\partial t}=P-\widehat{H}_{2}
\end{equation}
for $F$ and $\widehat{H}_{2}$ when $P$ is given. Equation \eqref{add1} implies that
\begin{equation}\label{add2}
-\omega \frac{\partial F}{\partial \theta}+\frac{\partial F}{\partial t}=P-\widehat{H}_{2}.
\end{equation}
Expand $F,P$ into Fourier series
\begin{equation}\label{add+1}
 F=\sum_{k\in\mathbb{Z}}F_{k}e^{ik\theta},  \ \
P=\sum_{k\in\mathbb{Z}}P_{k}e^{ik\theta},
\end{equation}
where the coefficients depend on $t$ and $r$.

Now, if we choose $\widehat{H}_{2}=[P]$, where $[P]$ stands for the mean value of $P$ with respect to the variable $\theta$, it is easy to show equation \eqref{add2} can always be solved. In fact, from \eqref{add2} and \eqref{add+1}, it follows that
$$ -ik\omega F_{k}+\frac{\partial F_{k}}{\partial t}=P_{k}.$$
Solving the above equation, one can obtain the $2\pi$-periodic solution
\begin{equation}\label{3.2}
F=\sum_{0\neq k\in\mathbb{Z}}e^{ik\omega t}\left(\int^{t}_{0}P_{k}e^{-ik\omega t}dt+\frac{1}{1-e^{2\pi ik\omega}}\int^{0}_{-2\pi}P_{k}e^{-ik\omega t}dt\right)e^{ik\theta},
\end{equation}
where we normalize  $F$ so that the mean value $[F]=0$.
Thus, $P-[P] \in \mbox{Range}(\mathcal{F})$, where $\mathcal{F}$ is the operator defined in Section 2, that is, $[P] \in \ker(\mathcal{F})$. From the hypotheses of Theorem 1.1, it follows that $\ker(\mathcal{F})=\mbox{span}\{H_{2}\},$
which yields that $[P]=0$.

Therefore, the Lie equation now is
\begin{equation} \label{outline1}
\{F,H_{2}\}+\frac{\partial F}{\partial t}=P,
\end{equation}
which can be completely solved. Hence, the new Hamiltonian is
\begin{equation}\label{Newequ}
H_+(r,\theta,t)=H_{2}+P_+,
\end{equation}
where
$$
P_{+}=\int^{1}_{0}\left\{(1-\varepsilon)\{H_{2},F\}+P-\frac{\partial F}{\partial t},F\right\} \circ X_{F}^{\varepsilon}d\varepsilon.
$$

If $P=O_{s}(r)$, it follows with \eqref{outline1} that $F=O_{s}(r),$  which implies that $P_{+}= O_{2s-1}(r),$ where $O_{s}(r)$ is a power series in $\sqrt{r}$ beginning with terms of degree at least $2s$, that is, $O_{s}(r)=O(r^s)$ as $r\to 0^+$.

From the above discussions, it is obvious that if we truncate $P$ with
\begin{equation}\label{outline2}
R=\sum_{i=2s}^{2(2s-1)}H_{i},
\end{equation}
we also have
\begin{eqnarray*}
P_{+}&=&\int^{1}_{0}\left\{(1-\varepsilon)\{H_{2},F\}+R-\frac{\partial F}{\partial t},F\right\} \circ X_{F}^{\varepsilon}d\varepsilon
+(P-R)\circ X_{F}^{1}
\\
&=& O_{2s-1}(r).
\end{eqnarray*}

Along this way, the initial Hamiltonian \eqref{intro3}
$$H=H_{2}+P=H_{2}+O_{s}$$
is transformed into
 $$H_{+}=H_{2}+P_{+}=H_{2}+O_{2s-1}.$$
Then we are in a position to carry out the following iteration process:
\begin{equation} \label{outline3}
\left\{
\begin{array}{ll}
s_{0}=2^{q}+1, s_{\nu}=2s_{\nu-1}-1=2^{q+\nu}+1, \nu=1,2,\cdots, \\[0.2cm]
H^{\nu}=H_{2}+P_{\nu}, \ P_{0}=O_{s_{0}}, \ P_{\nu}=O_{s_{\nu}},\\[0.2cm]
\mbox{the} \  \mbox{truncation} \
R_{\nu}=\sum_{i=2s_{\nu}}^{2(2s_{\nu}-1)}H_{i},\\[0.2cm]
F_{\nu}  \mbox{\ is  the   solution  of} \ \{F_{\nu},H_{2}\}+\frac{\partial F_{\nu}}{\partial t}=R_{\nu},\\[0.2cm]
\varphi_{\nu}=X_{F_{\nu}}^{1}, \ \psi_{0}=id, \ \psi_{\nu}=\varphi_{0}\circ\varphi_{1}\circ\cdots\circ\varphi_{\nu-1},\\[0.2cm]
\mbox{then} \ H^{\nu+1}=H_{2}+P_{\nu+1}=H_{2}+O_{2s_{\nu}-1},\\[0.2cm]
\mbox{and the  new  perturbation is } \\[0.2cm]
P_{\nu+1}=\int^{1}_{0}\left\{(1-\varepsilon)\{H_{2},F_{\nu}\}+R_{\nu}-\frac{\partial F_{\nu}}{\partial t},F_{\nu}\right\} \circ X_{F_{\nu}}^{\varepsilon}d\varepsilon \\[0.4cm]
 \ \ \ \ \ \ \ \ \ \ +(P_{\nu}-R_{\nu})\circ X_{F_{\nu}}^{1}.
\end{array}
\right.
\end{equation}

Without loss of generality, we can suppose that $q$ is as large as we want, because this can  be achieved
by using the existing formal transformation which transforms the Hamiltonian \eqref{intro3} into $H=H_{2}$
according to Theorem 1.1 and cutting the terms of degree higher than $2s_{0}-1.$ Then we obtain a convergent transformation, namely a polynomial, which transforms  the Hamiltonian \eqref{intro3} into $H=H_{2}+O_{s_{0}}.$


\section{The KAM step}

Now we consider one step of KAM process. Suppose $R$ defined in \eqref{outline2} is  analytic on $D_{\rho,\gamma}$ with
$|R|_{\rho,\gamma}\leq 1$. This reduction can be carried out by means of the change of variables $x=\varepsilon \hat{x}, y=\varepsilon \hat{y} $ in \eqref{intro2} with $\varepsilon$ sufficiently small. Then \eqref{intro2} has to be replaced by
$$\varepsilon^{-2}H(\varepsilon\hat{x},\varepsilon\hat{y},t)=H_{2}(\hat{x},\hat{y},t)+\varepsilon H_{3}(\hat{x},\hat{y},t)+\cdots+\varepsilon^{j-2}H_{j}(\hat{x},\hat{y},t)+\cdots,$$
and $P=O_{s}(r)=\varepsilon^{s-2}O_{s}(\hat{r})$ with $\hat{r}=\frac{1}{2}(\hat{x}^{2}+\hat{y}^{2})=\varepsilon^{-2}\hat{r}$.
In this way the domain of analyticity of $P$ can be made arbitrarily large and simultaneously the sup-norm of $P$ can be made arbitrarily small in a fixed domain.

The first result for the estimates of $F$ and its differential is the following.

\begin{lemma}\label{kam1}
Let $R$ be the truncation function of $P$ defined as \eqref{outline2}. Assume that $F$ solves the Lie equation
$$\{F,H_{2}\}+\frac{\partial F}{\partial t}=R,$$
and that the hypotheses of Theorem 1.1 are satisfied.
Then we have
\begin{equation*}\label{kamlm2}
|F|_{\rho-\delta,\gamma-\sigma}\leq c_{1}\sigma^{-\tau-1}e^{-s\delta}\delta^{-1},
\end{equation*}
where $\rho,\gamma \leq 1, \ 0<\delta<\frac{\rho}{4}, \ 0<\sigma<\frac{\gamma}{4}$, and $c_1$ is an absolute constant.
In addition, we have the following estimates
$$|F_{r}|_{\rho-2\delta,\gamma-\sigma}\leq c_{1}\sigma^{-\tau-1}e^{-s\delta}\delta^{-2},$$
$$|F_{\theta}|_{\rho-\delta,\gamma-2\sigma}\leq c_{1}\sigma^{-\tau-2}e^{-s\delta}\delta^{-1}.$$
\end{lemma}

\Proof From \eqref{3.2}, we can write the Fourier-Taylor series of $F$  as below:
$$F=\sum_{0\neq k\in\mathbb{Z}}\sum_{s\leq \frac{n}{2}\leq 2s-1}e^{ik\omega t}\left(\int^{t}_{0}R_{kn}e^{-ik\omega t}dt+\frac{1}{1-e^{2\pi ik\omega}}\int^{0}_{-2\pi}R_{kn}e^{-ik\omega t}dt\right)e^{ik\theta}r^{\frac{n}{2}}.$$
Now, we will estimate $F$ on $D_{\rho-\delta,\gamma-\sigma}$. Since $\omega$ satisfies \eqref{intro6}, the denominator $1-e^{2\pi ik\omega}$ can be estimated by
\begin{equation}\label{addpf1}
|1-e^{2\pi ik\omega}|=2|\sin(k\omega\pi)|=2\sin(|k\omega-l|\pi)\geq 2\frac{2}{\pi}|k\omega-l|\pi\geq\frac{4\alpha}{|k|^{\tau}}.
\end{equation}
On the other hand, since $R$ is real analytic, the Fourier coefficients $R_{kn}$ decay exponentially. Indeed, we have
\begin{equation}\label{addpf2}
|R_{kn}|_{\rho,\gamma}\leq A e^{-|k|\gamma}.
\end{equation}
Thus
\begin{equation}\label{addpf3}
\Big|\int _{0}^{t}R_{kn}e^{-ik\omega t}dt\Big|_{\rho,\gamma}\leq A e^{-|k|\gamma}\frac{2}{|ik\omega|}
\leq 2A e^{-|k|\gamma}\frac{|k|^{\tau}}{\alpha}.
\end{equation}
From \eqref{addpf1}, \eqref{addpf2} and \eqref{addpf3}, it follows that

\begin{eqnarray*}
 |F|_{\rho-\delta,\gamma-\sigma}
 &\leq&\Bigg|\sum_{k\neq 0}\left(2e^{-|k|\gamma}\frac{|k|^{\tau}}{\alpha}+
 \frac{\pi}{2}e^{-|k|\gamma}\frac{|k|^{\tau}}{\alpha}\right)e^{|k|(\gamma-\sigma)}
 \sum_{2s\leq n\leq2(2s-1)}r^{\frac{n}{2}}\Bigg|_{\rho-\delta,\gamma-\sigma}\\
 &\leq&\frac{4}{\alpha}\sum_{k\neq 0}|k|^{\tau}e^{-|k|\gamma}e^{|k|(\gamma-\sigma)}
 \sum_{2s\leq n\leq2(2s-1)}(\varrho-\delta)^{\frac{n}{2}}\\
 &\leq&\frac{8}{\alpha}\sum_{k=1}^{\infty}k^{\tau}e^{-k\sigma}\left(1-\frac{\delta}{\rho}\right)^{s}
 \sum_{0\leq n\leq 2s-2}\left(1-\frac{\delta}{\rho}\right)^{\frac{n}{2}}\\
 &\leq&\frac{8}{\alpha}\sum_{k=1}^{\infty}k^{\tau}e^{-k\sigma}e^{-s\frac{\delta}{\rho}}
 \left(\frac{\delta}{\rho}\right)^{-1}\\
 &\leq&\frac{8}{\alpha}\sum_{k=1}^{\infty}k^{\tau}e^{-k\sigma}e^{-s\delta}\delta^{-1}\\
 &\leq& c_{1}\sigma^{-\tau-1}e^{-s\delta}\delta^{-1}.
\end{eqnarray*}
Moreover, with Cauchy's estimate we get
$$|F_{r}|_{\rho-2\delta,\gamma-\sigma}\leq c_{1}\sigma^{-\tau-1}e^{-s\delta}\delta^{-2},$$
$$|F_{\theta}|_{\rho-\delta,\gamma-2\sigma}\leq c_{1}\sigma^{-\tau-2}e^{-s\delta}\delta^{-1}.$$
\qed

The following lemma will be important in this work.

\begin{lemma}\label{kamth}
Under the conditions of Lemma 4.1, if
\begin{equation}\label{th1}
c_{1}\sigma^{-\tau-2}e^{-s\delta}\delta^{-2}\leq 1,
\end{equation}
then there exists a symplectic, $2\pi$-periodic transformation
$\varphi: D_{\rho-3\delta,\gamma-3\sigma}\rightarrow D_{\rho,\gamma},$ which transforms $H=H_{2}+P=H_{2}+O_{s}(r)$ into
\begin{equation*}
H_{+}=H_{2}+P_{+}=H_{2}+O_{2s-1}(r).
\end{equation*}
Moreover, let us denote $\varphi=id+\Delta\varphi, \Delta\psi=\psi\circ\varphi-\psi,$ if there exists a consistent bound for $\psi'$, i.e., there exists a constant $B$ such that $|\psi'|_{\rho,\gamma}\leq B,$ then
\begin{equation}\label{kam4}
\big|\Delta\psi\big|_{\rho-3\delta,\gamma-3\sigma}
 \leq Bc_{1}\sigma^{-\tau-2}e^{-s\delta}\delta^{-2},
\end{equation}
\begin{equation}\label{kam5}
\big|\Delta\psi'\big|_{\rho-4\delta,\gamma-4\sigma}
 \leq Bc_{1}\sigma^{-\tau-3}e^{-s\delta}\delta^{-3}.
\end{equation}
\end{lemma}
\Proof By Lemma 4.1, we know
$$|F_{r}|_{\rho-2\delta,\gamma-\sigma}\leq c_{1}\sigma^{-\tau-1}e^{-s\delta}\delta^{-2},$$
$$|F_{\theta}|_{\rho-\delta,\gamma-2\sigma}\leq c_{1}\sigma^{-\tau-2}e^{-s\delta}\delta^{-1}.$$
From the assumption \eqref{th1}, it follows that
$$|F_{r}|_{\rho-2\delta,\gamma-\sigma}\leq c_{1}\sigma^{-\tau-1}e^{-s\delta}\delta^{-2}\leq \sigma,$$
$$|F_{\theta}|_{\rho-\delta,\gamma-2\sigma}\leq c_{1}\sigma^{-\tau-2}e^{-s\delta}\delta^{-1}\leq \delta.$$
Thus $|F_{r}|\leq \sigma$ and $|F_{\theta}|\leq \delta$ on $D_{\rho-2\delta,\gamma-2\sigma}$.

Now we define $\varphi$ as the real analytic time-1-map of the flow $X_{F}^{\varepsilon}$, with equations of motion
$$\dot{r}=F_{\theta}, \ \dot{\theta}=-F_{r}.$$
Therefore, $\varphi=X_{F}^{1}$ is well defined on $D_{\rho-3\delta,\gamma-3\sigma}$, with
$$\varphi=X_{F}^{1}:D_{\rho-3\delta,\gamma-3\sigma}\rightarrow D_{\rho-2\delta,\gamma-2\sigma}\subset D_{\rho,\gamma}.$$
With the discussion in Section 3,  we know that the Hamiltonian $H=H_{2}+P=H_{2}+O_{s}(r)$ is transformed into
$H_{+}=H_{2}+P_{+}=H_{2}+O_{2s-1}(r).$

Now we estimate $|\Delta\psi|$ and $|\Delta\psi'|$.
By Lemma 4.1, we have
\begin{eqnarray*}
|\Delta\varphi|_{\rho-3\delta,\gamma-3\sigma}&=&|\varphi-id|_{\rho-3\delta,\gamma-3\sigma}\\[0.2cm]
 &\leq&\max\left\{|F_{\theta}|_{\rho-\delta,\gamma-2\sigma},|F_{r}|_{\rho-2\delta,\gamma-\sigma}\right\}\\[0.2cm]
 &\leq& c_{1}\sigma^{-\tau-2}e^{-s\delta}\delta^{-2},
\end{eqnarray*}
and
\begin{eqnarray*}
|\Delta\psi|_{\rho-3\delta,\gamma-3\sigma}&=&|\psi\circ\varphi-\psi|_{\rho-3\delta,\gamma-3\sigma}\\[0.32cm]
&=&\Big|\int_{0}^{1}\frac{d}{d\varepsilon}\psi(id+\varepsilon\Delta\varphi)d\varepsilon\Big|_{\rho-3\delta,\gamma-3\sigma}\\[0.32cm]
 &\leq&|\psi'|_{\rho-3\delta,\gamma-3\sigma} \ |\Delta\varphi|_{\rho-3\delta,\gamma-3\sigma}\\[0.32cm]
 &\leq& Bc_{1}\sigma^{-\tau-2}e^{-s\delta}\delta^{-2}.
\end{eqnarray*}
With Cauchy's estimate we get
$$|\Delta\psi'|_{\rho-4\delta,\gamma-4\sigma}\leq Bc_{1}\sigma^{-\tau-3}e^{-s\delta}\delta^{-3}.$$
\qed


\section{The iteration process and the proof of Theorem 1.1}

In this Section, we will give the proof of Theorem 1.1. With the discussion in Section 4, the iteration process \eqref{outline3} now runs as follow:
\begin{equation} \label{pf1}
\left\{
\begin{array}{ll}
H^{\nu}=H_{2}+P_{\nu}=H_{2}+O_{s_{\nu}}, \ s_{\nu}=2^{q+\nu}+1,\\[0.32cm]
\rho_{0}=1, \ \rho_{\nu+1}=\rho_{\nu}-4\delta_{\nu}, \ \gamma_{0}=1, \ \gamma_{\nu+1}=\gamma_{\nu}-4\sigma_{\nu},\\[0.32cm]
\epsilon_{\nu}=Bc_{1}\sigma_{\nu}^{-\tau-3}e^{-s_{\nu}\delta_{\nu}}\delta_{\nu}^{-3},\\[0.32cm]
|\psi_{\nu+1}-id|_{\rho_{\nu+1}\gamma_{\nu+1}}
\leq\sum_{\mu=0}^{\nu}\epsilon_{\mu}\leq\frac{1}{2}, \\[0.32cm]
|\psi'_{\nu+1}-E|_{\rho_{\nu+1},\gamma_{\nu+1}}\leq\sum_{\mu=0}^{\nu}\epsilon_{\mu}\leq\frac{1}{2}, \ \nu=0,1,\cdots,
\end{array}
\right.
\end{equation}
where $B=\frac{3}{2}$, $id$ is the identity map, $E$ is the identity matrix, and $\delta_{\nu},\sigma_{\nu},\epsilon_{\nu}$ are positive constants.

From \eqref{pf1}, it follows that
\begin{equation} \label{pf2}
\rho_{\nu}=\rho_{0}-4\sum_{\mu=0}^{\nu-1}\delta_{\mu},\ \ \ \ \gamma_{\nu}=\gamma_{0}-4\sum_{\mu=0}^{\nu-1}\sigma_{\mu}.
\end{equation}
In order to obtain a common domain of definition $D_{\rho_{\infty},\gamma_{\infty}}\subset D_{\rho_{\nu},\gamma_{\nu}}, \ \nu=0,1,\cdots,$ it is necessary that $\rho_{\nu}\geq \rho_{\infty}>0, \ \gamma_{\nu}\geq \gamma_{\infty}>0,$
or, say
\begin{equation} \label{pf3}
\sum_{\mu=0}^{\infty}\delta_{\mu}\leq \frac{1}{8},\ \ \sum_{\mu=0}^{\infty}\sigma_{\mu}\leq \frac{1}{8}.
\end{equation}
Moreover, according to \eqref{pf1}, the following inequality must hold
\begin{equation} \label{pf4}
\sum_{\mu=0}^{\infty}\epsilon_{\mu}\leq \frac{1}{2}.
\end{equation}

Now let us assume that we are in a position to determine sequences of positive numbers $\delta_{0}, \delta_{1}, \cdots,$ $\sigma_{0}, \sigma_{1}, \cdots$ and $\epsilon_{0}, \epsilon_{1}, \cdots $ in a way that \eqref{pf1} is valid for the process \eqref{outline3} and that \eqref{pf3} and \eqref{pf4} are satisfied. Then the transformation $\psi_{\nu} $ and its differential are uniformly  bound on $D_{\frac{1}{2},\frac{1}{2}}$. Therefore there exists a subsequence $\{\psi_{\nu_{i}}\}$ which uniformly converges to an analytic function $\psi$ on $D_{\frac{1}{2},\frac{1}{2}}$, i.e., $\psi_{\nu_{i}}\rightarrow \psi$. Obviously, $\psi$ is symplectic, $2\pi$-periodic.

Now we turn to the proof of the estimate in \eqref{pf1} by the complete induction. We will still use the notations in Section 4. Obviously, Lemma 4.2 can be applied. In fact, $c_{1}\sigma_{\nu}^{-\tau-2}e^{-s_{\nu}\delta_{\nu}}\delta_{\nu}^{-2}\leq\epsilon_{\nu}\leq\frac{1}{2}<1, \nu=0,1,\cdots.$
Thus by Lemma 4.2, $\psi_{\nu}$ given by  \eqref{outline3} are well defined on $ D_{\rho_{\nu},\gamma_{\nu}} , \nu=0,1,\cdots.$

If $\nu=0,$ then $\psi_{1}=\varphi_{0}, \  \varphi_{0}=id+\Delta\varphi_{0}$, with Lemma 4.2,
\begin{eqnarray*}
 |\psi_{1}-id|_{\rho_{1},\gamma_{1}}&=&\big|\Delta\varphi_{0}\big|_{\rho_{0}-4\delta_{0},\gamma_{0}-4\sigma_{0}}\leq
 \big|\Delta\varphi_{0}\big|_{\rho_{0}-3\delta_{0},s_{0}-3\sigma_{0}} \\[0.3cm]
 &\leq&c_{1}\sigma_{0}^{-\tau-2}e^{-s_{0}\delta_{0}}\delta_{0}^{-2}\leq \epsilon_{0}.
\end{eqnarray*}
and with Cauchy's estimate, we have
\begin{eqnarray*}
 |\psi'_{1}-E|_{\rho_{1},\gamma_{1}}&=&\big|\Delta\varphi'_{0}\big|_{\rho_{0}-4\delta_{0},\gamma_{0}-4\sigma_{0}}\\[0.3cm]
 &\leq&c_{1}\sigma_{0}^{-\tau-3}e^{-s_{0}\delta_{0}}\delta_{0}^{-3}\leq \epsilon_{0}.
\end{eqnarray*}
 Thus the beginning step of \eqref{pf1} holds. Suppose that the iteration process \eqref{pf1} holds for the $\nu$-th step, we consider the $(\nu+1)$-th step. Since
$$\Delta\psi_{\nu+1}=\psi_{\nu+1}\circ\varphi_{\nu+1}-\psi_{\nu+1},$$
we know
\begin{eqnarray*}
\psi_{\nu+2}-id &=& \psi_{\nu+1}\circ\varphi_{\nu+1}-id \\[0.3cm]
&=&\psi_{\nu+1}\circ\varphi_{\nu+1}-\psi_{\nu+1}+\psi_{\nu+1}-id \\[0.3cm]
 &=& \Delta\psi_{\nu+1}+(\psi_{\nu+1}-id).
 \end{eqnarray*}
Under the assumption in the induction process, we have
$$|\psi'_{\nu+1}-E|_{\rho_{\nu+1},\gamma_{\nu+1}}\leq\sum_{\mu=0}^{\nu}\epsilon_{\mu}\leq\frac{1}{2},$$
and
$$|\psi'_{\nu+1}|_{\rho_{\nu+1},\gamma_{\nu+1}}\leq \frac{1}{2}+|E|_{\rho_{\nu+1},\gamma_{\nu+1}}\leq\frac{3}{2}=B,$$
thus there exists a consistent bound $B=\frac{3}{2}$  for  $\psi'$ on $D_{\rho_{\nu+1},\gamma_{\nu+1}}$.
From \eqref{kam4} and \eqref{kam5} in Lemma 4.2, it follows that
 $$\big|\Delta\psi_{\nu+1}\big|_{\rho_{\nu+2},\gamma_{\nu+2}}, \ \big|\Delta\psi'_{\nu+1}\big|_{\rho_{\nu+2},\gamma_{\nu+2}}\leq \epsilon_{\nu+1}.$$
Thus we get
\begin{eqnarray*}
\big|\psi_{\nu+2}-id\big|_{\rho_{\nu+2},\gamma_{\nu+2}}&=&
\big|\Delta\psi_{\nu+1}+(\psi_{\nu+1}-id)\big|_{\rho_{\nu+2},\gamma_{\nu+2}}\\[0.3cm]
 &\leq&\big|\Delta\psi_{\nu+1}\big|_{\rho_{\nu+2},\gamma_{\nu+2}}+\big|\psi_{\nu+1}-id\big|_{\rho_{\nu+2},\gamma_{\nu+2}}\\[0.3cm]
 &\leq&\epsilon_{\nu+1}+\sum_{\mu=0}^{\nu}\epsilon_{\mu}=\sum_{\mu=0}^{\nu+1}\epsilon_{\mu},
 \end{eqnarray*}
and
\begin{eqnarray*}
\big|\psi'_{\nu+2}-E\big|_{\rho_{\nu+2},\gamma_{\nu+2}}&=&
\big|\Delta\psi'_{\nu+1}+(\psi'_{\nu+1}-E)\big|_{\rho_{\nu+2},\gamma_{\nu+2}}\\[0.3cm]
 &\leq&\big|\Delta\psi'_{\nu+1}\big|_{\rho_{\nu+2},\gamma_{\nu+2}}+\big|\psi'_{\nu+1}-E\big|_{\rho_{\nu+2},\gamma_{\nu+2}}\\[0.3cm]
 &\leq&\epsilon_{\nu+1}+\sum_{\mu=0}^{\nu}\epsilon_{\mu}=\sum_{\mu=0}^{\nu+1}\epsilon_{\mu}.
 \end{eqnarray*}

As a result, the estimates in \eqref{pf1} are valid for $\nu+1$ instead of $\nu$ and we are left with the problem of finding a sequence of
positive constants $\delta_{0}, \delta_{1},\cdots, \sigma_{0}, \sigma_{1}, \cdots$ such that the inequalities \eqref{pf3} and \eqref{pf4}
are satisfied.

Now we turn to find  appropriate parameters.
Choose
$$\sigma_{\nu}=\frac{1}{2^{\nu+4}}, \ \nu=0,1,\cdots ,$$
then$$\sum_{\nu=0}^{\infty}\sigma_{\nu}=\frac{1}{8}.$$
Thus, $\epsilon_{\nu}$ has the form
\begin{eqnarray*}
\epsilon_{\nu}&=&Bc_{1}\sigma_{\nu}^{-\tau-3}e^{-s_{\nu}\delta_{\nu}}\delta_{\nu}^{-3}
=Bc_{1}2^{(\tau+3)(\nu+4)}e^{-s_{\nu}\delta_{\nu}}\delta_{\nu}^{-3}\\[0.3cm]
&=&c_{2}2^{(\tau+3)(\nu+4)}e^{-s_{\nu}\delta_{\nu}}\delta_{\nu}^{-3},
\end{eqnarray*}
where $c_{2}=Bc_{1}.$ Without loss of generality, we may suppose $c_{2}\geq 1$. We get from \eqref{pf4}
the necessary condition
$$2^{(\tau+3)(\nu+4)}e^{-s_{\nu}\delta_{\nu}}\leq 1,$$
or
$$\delta_{\nu}\geq\frac{(\tau+3)(\nu+4)\ln2}{s_{\nu}}=\frac{(\tau+3)(\nu+4)\ln2}{2^{q+\nu}+1}.$$
Choose
$$\delta_{\nu}=\frac{4(\tau+3)(\nu+q)\ln2}{s_{\nu}}=\frac{4(\tau+3)(\nu+q)\ln2}{2^{q+\nu}+1}, \
\nu=0,1,\cdots ,$$
then
\begin{eqnarray*}
\sum_{\nu=0}^{\infty}\delta_{\nu}&=&\sum_{\nu=0}^{\infty}\frac{4(\tau+3)(\nu+q)\ln2}{2^{q+\nu}+1}
\leq\sum_{\nu=0}^{\infty}\frac{4(\tau+3)(\nu+q)}{2^{q+\nu}}\\
&=&\sum_{\nu=0}^{\infty}4(\tau+3)\frac{\nu+q}{2^{q+\nu}}\rightarrow 0 , \ \ \mbox{as}\ \ q\rightarrow\infty.
\end{eqnarray*}
With the  parameters $\delta_{\nu}, \ \sigma_{\nu}$, we have
\begin{eqnarray*}
\epsilon_{\nu}&=&Bc_{1}\sigma_{\nu}^{-\tau-3}e^{-s_{\nu}\delta_{\nu}}\delta_{\nu}^{-3}
= Bc_{1}2^{(\tau+3)(\nu+4)}e^{-s_{\nu}\delta_{\nu}}\delta_{\nu}^{-3}\\[0.3cm]
&\leq&c_{2}2^{(\tau+3)(\nu+q)}e^{-4(\tau+3)(\nu+q)\ln2}\left(\frac{4(\tau+3)(\nu+q)\ln2}{2^{q+\nu}+1}\right)^{-3}\\[0.3cm]
&=&c_{2}\left(2^{(\tau+3)(\nu+q)}\frac{4(\tau+3)(q+\nu)\ln2}{2^{q+\nu}+1}\right)^{-3}\\[0.3cm]
&\leq&c_{2}\left(4(\tau+3)\ln2\right)^{-3}(q+\nu)^{-3}\\[0.3cm]
&=&c_{2}\left(4(\tau+3)\ln2\right)^{-3}(q^{3}+\nu^{3}+3q^{2}\nu+3q\nu^{2})^{-1}\\[0.3cm]
&\leq &c_{3}\frac{1}{3q}\frac{1}{\nu^{2}},
\end{eqnarray*}
where we choose $q\geq 4, \ c_{3}=c_{2}\left(4(\tau+3)\ln2\right)^{-3}.$
Thus
$$\sum_{\nu=0}^{\infty}\epsilon_{\nu}\leq c_{3}\frac{1}{3q}\sum_{\nu=0}^{\infty}\frac{1}{\nu^{2}}\rightarrow 0 , \ \mbox{as}\
\ q\rightarrow\infty.$$

Therefore, the  parameters $\delta_{\nu}, \ \sigma_{\nu}$ chosen as above satisfy \eqref{pf3} and
\eqref{pf4} if $q$ is large enough. Up to now we have completed the proof of Theorem 1.1.

\section*{Acknowledgments} We would like to thank the referees for their valuable comments which have led to an improvement of the presentation of this paper.


\section*{References}
\bibliographystyle{elsarticle-num}

\begin{thebibliography}{99}

\bibitem{Arnold} V. I. Arnold,  {\em The stability of the equilibrium position of a Hamiltonian system of ordinary differential equations in the general elliptic case}, Soviet Math. Dokl., 2 (1961), 247-249.




\bibitem{Bardin} B. S. Bardin, V. Lanchares, {\em On the stability of periodic Hamiltonian systems with one degree of freedom in the case of degeneracy}, Regul. Chaotic Dyn., 20 (2015), 627-648.


\bibitem{Caraballo} T. Caraballo, P. E. Kloeden, {\em Stability of equilibria and fixed points of conservative systems}, Nonlinearity, 12 (1999), 1351-1362.

 \bibitem{Chu} J. Chu, M. Zhang, {\em  Rotation numbers and Lyapunov stability of elliptic periodic solutions}, Discrete Contin. Dyn. Syst., 21 (2008), 1071-1094.

 \bibitem{Lei} J. Lei, X. Li, P. Yan, M. Zhang, {\em Twist character of the least amplitude periodic solution of the forced pendulum}, SIAM J. Math. Anal., 35 (2003), 844-867.

\bibitem{Liu1} B. Liu, {\em The stability of the equilibrium of a conservative system}, J. Math. Anal. Appl., 202 (1996), 133-149.

\bibitem{Mansilla} J. E. Mansilla, C. Vidal, {\em Stability of equilibrium solutions in the critical case of even-order resonance in periodic Hamiltonian systems with one degree of freedom}, Celestial Mech. Dynam. Astronom., 116 (2013), 265-277.


\bibitem{Markeyev} A. P. Markeyev, {\em The critical case of fourth-order resonance in a Hamiltonian system with one degree of freedom}, J. Appl. Maths. Mechs., 61 (1997), 355-361.


\bibitem{Meyer}K. Meyer, G. R. Hall, {\em Introduction to Hamiltonian dynamical systems and the n-body problem}, New York. Springer-Verlag. Applied mathematical science, 1992.

 \bibitem{Moser0} J. Moser,  {\em New aspects in the theory of stability of Hamiltonian systems}, Comm. Pure Appl. Math., 11 (1958), 81-114.

\bibitem{Moser} J. Moser, {\em On invariant curves of area preserving mappings of an annulus},
Nachr. Akad. Wiss. Gott., Math. Phys. Kl., (1962), 1-20.

\bibitem{Moser1} J. Moser, {\em Stable and Random Motions in Dynamical Systems},
 University Press, Princeton, 2001.


\bibitem{Russmann1} H. R\"{u}ssmann,  {\em Stability of elliptic fixed points of analytic area-preserving mappings under the Bruno condition},
Ergod. Th. \& Dynam. Sys., 22 (2002), 1551-1573.

\bibitem{Ortega0} R. Ortega, {\em The twist coefficient of periodic solutions of a time-dependent Newton's equation}, J. Dynam. Differential Equations, 4 (1992), 651-665.

\bibitem{Ortega1} R. Ortega, {\em The stability of the equilibrium of a nonlinear Hill's equation}, SIAM J. Math. Anal., 25 (1994), 1393-1401.

\bibitem{Ortega2} R. Ortega, {\em Periodic solutions of a Newtonian equation: stability by the third approximation}, J. Differential Equations, 128 (1996), 491-518.

\bibitem{Ortega3} D. N\'{u}\~{e}z, R. Ortega, {\em Parabolic fixed points and stability criteria for nonlinear Hill's equation}, Z. Angew. Math. Phys., 51 (2000), 890-911.
 \bibitem{Ortega4} R. Ortega, {\em The stability of the equilibrium: a search for the right approximation}, Ten mathematical essays on approximation in analysis and topology, Elsevier B. V., Amsterdam, (2005), 215-234.

 \bibitem{Raissy1} J. Raissy, {\em Linearization of holomorphic germs with quasi-Brjuno fixed points}, Math. Z., 264 (2010), 881-900.

\bibitem{Raissy2} J. Raissy, {\em Brjuno conditions for linearization in presence of resonances}, Asymptotics in dynamics, geometry and PDEs; generalized Borel summation, 201-218, CRM Series, 12, Ed. Norm., Pisa, 2011.

  \bibitem{Raissy3} J. Raissy, {\em Holomorphic linearization of commuting germs of holomorphic maps}, J. Geom. Anal., 23 (2013), 1993-2019.


 \bibitem{Santos1} F. dos Santos, J. E. Mansilla, C. Vidal, {\em Stability of equilibrium solutions of autonomous and periodic Hamiltonian systems with n-degrees of freedom in the case of single resonance}, J. Dyn. Diff. Equat., 22 (2010), 805-821.


\bibitem{Santos2} F. dos Santos,  C. Vidal, {\em Stability of equilibrium solutions of autonomous and periodic Hamiltonian systems  in the case of multiple resonances}, J. Differential Equations, 258 (2015), 3880-3901.


\bibitem{Santos3} C. Vidal, F. dos Santos, {\em Stability of equilibrium solutions of  periodic Hamiltonian systems  under third and fourth order resonances}, Regul. Chaotic Dyn., 10 (2005), 95-110.



\bibitem{Siegel} C. L. Siegel,  J. K. Moser, {\em Lectures on celestial mechanics}, Springer, 1971.

\bibitem{Torres1} P. J. Torres, {\em Twist solutions of a Hill's equation with singular term}, Adv. Nonlinear Stud., 2 (2002), 279-287.

 \bibitem{Torres2} P. J. Torres, {\em Existence and stability of periodic solutions for second-order semilinear differential equations with a singular nonlinearity}, Proc. Roy. Soc. Edinburgh Sect., 137 (2007), 195-201.

\bibitem{Wu} Y. Wu, {\em Linearizability of quasi-periodic Hamiltonian systems in resonant case}, Adv. Math., (China) 41 (2012), 605-614.




















\end{thebibliography}

\end{document}